Analysis.

# About polynomials related to a quadratic equation.

## Roland Groux

*Lycée Polyvalent Rouvière, rue Sainte Claire Deville,
BP 1205. 83070 Toulon .Cedex. France.*

Email : roland.groux@orange.fr

**Abstract.**

We consider here a particular quadratic equation linking two elements of a C-Algebra. By analysing powers of the unknowns, it appears a double sequence of polynomials related to classical Bernoulli polynomials. We get the generating functions, integral forms and explicit formulas for the coefficients involving cosecant and tangent numbers. We also study the use of these polynomials for the calculation of some integral transforms.

## 1. Introduction and notations.

In what follows we consider the equation (1.1) : $\boxed{2ab = a^2 + b^2 + 1}$, with $a$, $b$ elements of any Algebra over the field of complex numbers.

We first prove by induction the general result:

(1.2) For every $n$ integer there exists two polynomials noted $A_n$ and $C_n$, satisfying the equalities : $\boxed{\begin{cases} ab^n = A_n(a) + C_n(b) \\ a^n b = C_n(a) + A_n(b) \end{cases}}$, with the leading term of $C_n(X)$ equal to $\dfrac{n}{n+1} X^{n+1}$.

For $n = 0$ we have $A_0(X) = X$ and $C_0(X) = 0$.

For $n = 1$ we have $A_1(X) = \dfrac{X^2}{2}$ and $C_1(X) = \dfrac{X^2 + 1}{2}$.

So we can write : $(R_n)$ : $ab^n = A_n(a) + C_n(b)$.

Multiplying on the right the terms of (1.1) by $b^n$ we get : $2ab^{n+1} = a^2 b^n + (b^2 + 1)b^n$.

Thanks to $(R_n)$ we deduce: $a^2 b^n = aA_n(a) + aC_n(b)$.

By combining the two previous equalities: (F) $2ab^{n+1} = aA_n(a) + aC_n(b) + (b^2+1)b^n$

Using the recurrence hypothesis: $C_n(b) = \sum_{k=0}^{k=n} \lambda_n^k b^k + \dfrac{n}{n+1} b^{n+1}$, with $ab^k = A_k(a) + C_k(b)$

So : $aC_n(b) = \sum_{k=0}^{k=n} \lambda_n^k (A_k(a) + C_k(b)) + \dfrac{n}{n+1} ab^{n+1}$ and the formula (F) above translate to:

$2ab^{n+1} = aA_n(a) + \dfrac{n}{n+1} ab^{n+1} + \sum_{k=0}^{k=n} \lambda_n^k (A_k(a) + C_k(b)) + (b^2+1)b^n$.

By reducing the terms we finally obtain :

$$\left(\frac{n+2}{n+1}\right)ab^{n+1} = aA_n(a) + \sum_{k=0}^{k=n}\lambda_n^k(A_k(a) + C_k(b)) + (b^2+1)b^n$$

By introducing the polynomials defined by :

$$(1.3) \quad \boxed{\begin{aligned} A_{n+1}(X) &= \left(\frac{n+1}{n+2}\right)\left(XA_n(X) + \sum_{k=0}^{k=n}\lambda_n^k A_k(X)\right) \\ C_{n+1}(X) &= \left(\frac{n+1}{n+2}\right)\left((X^2+1)X^n + \sum_{k=0}^{k=n}\lambda_n^k C_k(X)\right) \end{aligned}}, \text{ we get :}$$

$ab^{n+1} = A_{n+1}(a) + C_{n+1}(b)$ and the leading term of $C_{n+1}(X)$ is $\left(\frac{n+1}{n+2}\right)X^{n+2}$.

So the first equality of (1.2) is proved by induction. A similar approach easily sets to the second part of (1.2) with the same couple of polynomials.

Note that in terms of umbral composition [see 1], we can write (1.3) as :

$$A_{n+1}(X) = \left(\frac{n+1}{n+2}\right)(XA_n(X) + C_n(A(X)))$$

$$C_{n+1}(X) = \left(\frac{n+1}{n+2}\right)((X^2+1)X^n + C_n(C(X)))$$

Here the first terms :

$$A_0(X) = X; \quad A_1(X) = \frac{X^2}{2}; \quad A_2(X) = \frac{X^3 + X}{3}; \quad A_3(X) = \frac{X^4 + 2X^2}{4}; \quad A_4(X) = \frac{3X^5 + 10X^3 + 7X}{15}$$

$$A_5(X) = \frac{X^6 + 5X^4 + 7X^2}{6}; \quad A_6(X) = \frac{3X^7 + 21X^5 + 49X^3 + 31X}{21}$$

$$A_7(X) = \frac{X^8}{8} + \frac{7}{6}X^6 + \frac{49}{12}X^4 + \frac{31}{6}X^2; \quad A_8(X) = \frac{X^9}{9} + \frac{4}{3}X^7 + \frac{98}{15}X^5 + \frac{124}{9}X^3 + \frac{127}{15}X$$

$$C_0(X) = 0; \quad C_1(X) = \frac{X^2+1}{2}; \quad C_2(X) = \frac{2X^3 + 2X}{3}; \quad C_3(X) = \frac{3X^4 + 4X^2 + 1}{4}$$

$$C_4(X) = \frac{4}{5}X^5 + \frac{4}{3}X^3 + \frac{8}{15}X; \quad C_5(X) = \frac{5}{6}X^6 + \frac{5}{3}X^4 + \frac{4}{3}X^2 + \frac{1}{2}$$

$$C_6(X) = \frac{6}{7}X^7 + 2X^5 + \frac{8}{3}X^3 + \frac{32}{21}X; \quad C_7(X) = \frac{7}{8}X^8 + \frac{7}{3}X^6 + \frac{14}{3}X^4 + \frac{16}{3}X^2 + \frac{17}{8}$$

## 2. Basic relations.

2.1 An elementary example.

Over the field of complex numbers we consider: $a = x - i$; $b = x$.

The general equality $ab^n = A_n(a) + C_n(b)$ here leads to: $(x-i)x^n = A_n(x-i) + C_n(x)$, for every $x$ complex.

So we deduce a first elementary equality linking $A_n$ and $C_n$ :

$$\boxed{A_n(X-i) + C_n(X) = (X-i)X^n} \quad (2.1.1)$$

Now, with $a = x+i$; $b = x$, we get: $\boxed{A_n(X+i) + C_n(X) = (X+i)X^n}$ (2.1.2)

By difference : $A_n(X+i) - A_n(X-i) = 2iX^n$ (2.1.3)

Note then $B_n(X)$ the bernoulli polynomial of order $n$ and recall the well know identity :

$B_n(x+1) - B_n(x) = nx^{n-1}$. [see 2]

From this it is easy to prove that $U_n(X) = \dfrac{(2i)^{n+1}}{n+1} B_{n+1}(\dfrac{X+i}{2i})$ satisfies the difference equation

$U_n(X+i) - U_n(X-i) = 2iX^n$, similar to (2.1.3). So the difference between $A_n(X)$ and $U_n(X)$ is constant.

For the initial value, it is obvious thanks the formula (1.3) that $A_n(0) = 0$ for every $n$.

We conclude the egality :(2.1.3) $\boxed{A_n(X) = \dfrac{(2i)^{n+1}}{n+1}[B_{n+1}(\dfrac{X+i}{2i}) - B_{n+1}(\dfrac{1}{2})]}$

We can perform a similar approach with the set of $C_n(X)$.

From (2.1.1) and (2.1.2) we get : $\begin{cases} A_n(X) + C_n(X+i) = X(X+i)^n \\ A_n(X) + C_n(X-i) = X(X-i)^n \end{cases}$

By difference : (2.1.4) $C_n(X+i) - C_n(X-i) = X[(X+i)^n - (X-i)^n]$

The polynomial $V_n(X) = X^n(X-i) - \dfrac{(2i)^{n+1}}{n+1} B_{n+1}(\dfrac{X}{2i})$ satisfies also to the same difference equation. Or for every $n$ integer we deduce from (1.3) the value $C_n(i) = 0$.

So we conclude the identity (2.1.5) $\boxed{C_n(X) = \dfrac{(2i)^{n+1}}{n+1}[-B_{n+1}(\dfrac{X}{2i}) + B_{n+1}(\dfrac{1}{2})] + X^n(X-i)}$

2.2 Links with the Euler polynomials. [see 2]

Recall the well known identity linking $E_n(X)$ the Euler polynomial of order $n$ with the Bernoulli polynomial $B_n(X)$ : $E_n(X) = \dfrac{2^{n+1}}{n+1}[B_{n+1}(\dfrac{X}{2} + \dfrac{1}{2}) - B_n(\dfrac{X}{2})]$. Using this formula through the addition of (2.1.3) and (2.1.4) we easily obtain :

$$(2.2.1) \boxed{E_n(X) = X^n - X^{n+1} + (-i)^{n+1}[A_n(ix) + C_n(ix)]}$$

## 3. Explicitations of the coefficients.

We recall that a generative function of a sequence of polynomials $W_n(x)$ is defined by the formal series: $h(t,x) = \sum_{n=0}^{n=\infty} \frac{W_n(x)}{n!} t^n$. [see 3]

In the case of Bernoulli polynomials it is well known that $h(t,x) = \frac{te^{xt}}{e^t - 1}$.

So, from (2.1.3) it is easy to deduce the generative function of the sequence $A_n(X)$.

We get: $\boxed{f(t,x) = \sum_{n=0}^{n=\infty} \frac{A_n(x)}{n!} t^n = \frac{e^{xt} - 1}{\sin(t)}}$. (3.1)

> `f:=(exp(x*t)-1)/sin(t) :`

> `series(f,t=0,8);`

$$x + \frac{x^2}{2}t + \left(\frac{1}{6}x^3 + \frac{1}{6}x\right)t^2 + \left(\frac{1}{24}x^4 + \frac{1}{12}x^2\right)t^3 + \left(\frac{1}{120}x^5 + \frac{7}{360}x + \frac{1}{36}x^3\right)t^4 +$$
$$\left(\frac{1}{720}x^6 + \frac{7}{720}x^2 + \frac{1}{144}x^4\right)t^5 + \left(\frac{1}{5040}x^7 + \frac{31}{15120}x + \frac{7}{2160}x^3 + \frac{1}{720}x^5\right)t^6 + O(t^7)$$

Thanks to this result we can make explicit the coefficients of $A_n(X)$. We recall the definition of cosecant numbers: [see 4] $\frac{x}{\sin(x)} = \sum_{n=0}^{n=\infty} \frac{cs(n)}{n!} x^n$. They can be expressed using the Bernoulli numbers $\beta_n = B_n(0)$ by: $cs(n) = (-1)^{\frac{n}{2}+1} (2^n - 2)\beta_n$

Then we get: $A_n(x) = \frac{1}{n+1} \sum_{k=1}^{k=n+1} C_{n+1}^k cs(n+1-k) x^k$ ; $\alpha_n^k = \frac{1}{n+1} C_{n+1}^k cs(n+1-k)$ (3.2)

In the same way we deduce from (2.1.4) $\boxed{g(x,t) = \sum_{n=0}^{n=\infty} \frac{C_n(x)}{n!} t^n = xe^{tx} + \frac{1 - e^{tx}\cos(t)}{\sin(t)}}$. (3.3)

> `g:=x*exp(t*x)+(1-exp(t*x)*cos(t))/sin(t) :`
> `series(g,t=0,8);`

$$\left(\frac{x^2}{2} + \frac{1}{2}\right)t + \left(\frac{1}{3}x^3 + \frac{1}{3}x\right)t^2 + \left(\frac{1}{8}x^4 + \frac{1}{24} + \frac{1}{6}x^2\right)t^3 + \left(\frac{1}{30}x^5 + \frac{1}{45}x + \frac{1}{18}x^3\right)t^4 +$$
$$\left(\frac{1}{144}x^6 + \frac{1}{72}x^4 + \frac{1}{240} + \frac{1}{90}x^2\right)t^5 + \left(\frac{1}{360}x^5 + \frac{1}{840}x^7 + \frac{1}{270}x^3 + \frac{2}{945}x\right)t^6 + O(t^7)$$

If we consider the series: $\tan(\frac{x}{2}) = \sum_{n=1}^{n=\infty} \frac{d_n}{n!} x^n$, with $d_n = 2(-1)^{\frac{n-1}{2}} (2^{n+1} - 1) \frac{\beta_{n+1}}{n+1}$

We deduce by addition of the two generative functions (3.1) ; (3.3) :

$$A_n(x) + C_n(x) = x^{n+1} + \sum_{k=0}^{k=n-1} C_n^k d_{n-k} x^k$$

And finally, using the formula (3.2) we get the coefficient $\lambda_n^k$ of $C_n(X)$ :

(3.4) For $k \notin \{0, n, n+1\}$ $\lambda_n^k = (-1)^{\frac{n-k-1}{2}} C_n^k 2^{n-k+1} \frac{\beta_{n-k+1}}{n-k+1}$ ;

$\lambda_n^0 = 2(2^{n+1}-1)(-1)^{\frac{n-1}{2}} \frac{\beta_{n+1}}{n+1}$ ; $\lambda_n^n = 0$ ; $\lambda_n^{n+1} = \frac{n}{n+1}$ ;

## 4. Applications to integral calculus.

4.1 Introduction.

We consider here the classic Hilbert space $L^2([0,1])$. We note $\varphi_0$ the function $x \mapsto \varphi_0(x) = \ln(\frac{x}{1-x})$ and $T$ the operator defined by $f(x) \mapsto g(x) = \int_0^1 \frac{f(t)-f(x)}{t-x}$.

Recall the following result: If $f$, $\varphi_0 \times f$ and $T(f)$ are both elements of $L^2([0,1])$, then we have the identity (4.1.1) : $T(2\varphi_0 \times f - T(f)) = (\varphi_0^2 + \pi^2) \times f$ [see 5]

If we introduce the operators $a = \frac{T}{\pi}$ and $f \mapsto b(f) = \frac{\varphi_0}{\pi} \times f$, then (4.1.1) leads to the quadratic equation (1.1) : $\boxed{2ab = a^2 + b^2 + 1}$

So we can consider in this case some consequences of the previous general results.

In what follows we assume that the assumptions for applying (4.1.1) remain valid during the successive iterations of the operators.

4.2 First consequences of $ab^n = A_n(a) + C_n(b)$.

The above formula gives here directly:

$$\boxed{T(\varphi_0^n f) = \sum_{k=0}^{k=n+1} \alpha_n^k \pi^{n+1-k} T^k(f) + (\sum_{k=0}^{k=n+1} \lambda_n^k \pi^{n+1-k} \varphi_0^k) \times f}$$ (4.2.1)

For $f = 1$ this leads to $T(\varphi_0^n) = \sum_{k=0}^{k=n+1} \lambda_n^k \pi^{n+1-k} \varphi_0^k$ , what is written as :

$$\int_0^1 \frac{\ln^n(\frac{u}{1-u}) - \ln^n(\frac{x}{1-x})}{u-x} du = \pi^{n+1} C_n(\frac{\varphi_0(x)}{\pi}).$$

By change of variables: $z = \ln(\frac{x}{1-x})$ and $t = \frac{u}{1-u}$, we obtain an integral form for the polynomial $C_n$, all over $\mathbb{R}$ : $\boxed{(e^z + 1)\int_0^{+\infty} \frac{\ln^n(t) - z^n}{(1+t)(t-e^z)} dt = \pi^{n+1} C_n(\frac{z}{\pi})}$. (4.2.2)

For $z = 0$ we have : $2\int_0^{+\infty} \frac{\ln^n(t)}{t^2-1} dt = \pi^{n+1} \lambda_n^0$.

Thanks to the expression of $\lambda_n^0 = 2(2^{n+1}-1)(-1)^{\frac{n-1}{2}}\frac{\beta_{n+1}}{n+1}$, and the transformation

$2\int_0^{+\infty}\frac{\ln^n(t)dx}{t^2-1} = 2\int_0^1(1+(-1)^{n+1})\frac{\ln^n(t)dx}{t^2-1}$, we get for $n$ odd integer the classical integral :

$$4\int_0^1\frac{\ln^{2n-1}(x)dx}{x^2-1} = \frac{(4^n-1)(-1)^{n-1}\beta_{2n}}{n}\pi^{2n}$$

### 4.3 Evaluations of integrals of type $\int_0^1 \varphi_0^n(x)f(x)dx$.

We show easily the followings formulas : (4.3.1)

$$T(xf(x)) = xT(f(x)) + \int_0^1 f(x)dx \ ; \ \ldots\ ; \ T^n(xf(x)) = xT^n(f(x)) + \int_0^1 T^{n-1}(f(x))dx$$

So we can evaluate $T(\varphi_0^n(x)xf(x))$ by two different ways :

- With (4.2.1) : $T(\varphi_0^n(x)xf(x)) = \sum_{k=1}^{k=n+1}\alpha_n^k\pi^{n+1-k}T^k(xf(x)) + (\sum_{k=0}^{k=n+1}\lambda_n^k\pi^{n+1-k}\varphi_0^k(x))\times xf(x)$

- With (4.3.1) : $T(x\varphi_0^n(x)f(x)) = xT(\varphi_0^n(x)f(x)) + \int_0^1 \varphi_0^n(x)f(x)dx$

From (4.2.1) : $T(\varphi_0^n(x)f(x)) = \sum_{k=1}^{k=n+1}\alpha_n^k\pi^{n+1-k}T^k(f(x)) + (\sum_{k=0}^{k=n+1}\lambda_n^k\pi^{n+1-k}\varphi_0^k(x))\times f(x)$

From (4.3.1) and every $k$ element of $\{1,\ldots,n-1\}$ : $T^k(xf(x)) = xT^k(f(x)) + \int_0^1 T^{k-1}(f(x))dx$.

So we can translate the two methods above as :

(1) : $T(\varphi_0^n(x)xf(x)) = \sum_{k=1}^{k=n+1}\alpha_n^k\pi^{n+1-k}[xT^k(f(x)) + \int_0^1 T^{k-1}(f(x))dx] + (\sum_{k=0}^{k=n+1}\lambda_n^k\pi^{n+1-k}\varphi_0^k(x))xf(x)$

(2) $T(\varphi_0^n(x)xf(x)) = \sum_{k=1}^{k=n+1}\alpha_n^k\pi^{n+1-k}xT^k(f(x)) + (\sum_{k=0}^{k=n+1}\lambda_n^k\pi^{n+1-k}\varphi_0^k(x))xf(x) + \int_0^1 \varphi_0^n(x)f(x)dx$

By comparing we deduce : $\boxed{\int_0^1 \varphi_0^n(x)f(x)dx = \sum_{k=1}^{k=n+1}\alpha_n^k\pi^{n+1-k}\int_0^1 T^{k-1}(f(x))dx}$ (4.3.2)

Apply this result to the function : $x \mapsto f(x) = \frac{1}{x+a}$

It is obvious that $f$ is an eigenvector of $T$, with the eigenvalue $\gamma_a = \ln(\frac{a}{1+a})$

So, for every $k \geq 1$ : $T^k(f(x)) = (\gamma_a)^k f(x)$ and we deduce :

$\int_0^1 T^{k-1}(f(x))dx = (\gamma_a)^{k-1} \int_0^1 \frac{dx}{x+a} = -(\gamma_a)^k$ . In this particular case (4.3.2) gives :

$$\int_0^1 \frac{\varphi_0^n(x)}{x+a} dx = -\sum_{k=1}^{k=n+1} \alpha_n^k \pi^{n+1-k} (\gamma_a)^k = -\pi^{n+1} A_n(\frac{\gamma_a}{\pi}) \quad (4.3.3)$$

By change of variable $t = \frac{x}{1-x}$, we get :

$\frac{1}{1+a} \int_0^{+\infty} \frac{\ln^n(t)dt}{(1+t)(t+\frac{a}{1+a})} = -\pi^{n+1} A_n(\frac{\gamma_a}{\pi})$. Now, with $z = \gamma_a = \ln(\frac{a}{1+a})$, we conclude an

integral form for the polynomial $A_n$ :

$$(1-e^z) \int_0^{+\infty} \frac{\ln^n(t)dt}{(t+1)(t+e^z)} = -\pi^{n+1} A_n(\frac{z}{\pi}) \quad (4.3.4)$$

By a classic method using residues, we can evaluate this integral with a recurrence.

After simplifications we get : $\sum_{k=0}^{k=n-1} C_n^k (2i)^{n-(k+1)} A_k(\frac{z}{\pi}) = (\frac{z}{\pi}+i)^n - (i)^n$

So the polynomials $A_n$ satisfy :

$$A_{n+1}(X) = \frac{1}{n+2}[(X+i)^{n+2} - (i)^{n+2} - \sum_{k=0}^{k=n} C_{n+2}^k (2i)^{n+1-k} A_k(X)] \quad (4.3.5)$$

4.4 Consequences of $a^n b = C_n(a) + A_n(b)$.

This is reflected here by :

$$T^n(\varphi_0 \times f) = \sum_{k=0}^{k=n+1} \lambda_n^k \pi^{n+1-k} T^k(f) + (\sum_{k=0}^{k=n+1} \alpha_n^k \pi^{n+1-k} \varphi_0^k) \times f \quad (4.4.1)$$

For $f = 1$ we get : $T^n(\varphi_0) = \lambda_n^0 \pi^{n+1} + \sum_{k=0}^{k=n+1} \alpha_n^k \pi^{n+1-k} \varphi_0^k \quad (4.4.2)$

Calculating as before $T^n(\varphi_0 \times f \times x)$ by two methods, we obtain a similar formula for the integrals over $[0,1]$ of the transforms $T^k(\varphi_0 \times f)$ :

$$\int_0^1 T^{n-1}(\varphi_0(x)f(x))dx = \sum_{k=1}^{k=n+1} \lambda_n^k \pi^{n+1-k} \int_0^1 T^{k-1}(f(x))dx \quad (4.4.3)$$

The particular case $f=1$, leads to : $\boxed{\int_0^1 T^{n-1}(\varphi_0(x))dx = \lambda_n^1 \pi^n}$ (4.4.5)

For odd values of the power : $\int_0^1 T^{2n-1}(\ln(\frac{x}{1-x}))dx = \lambda_{2n}^1 \pi^{2n} = (-1)^{n+1}(4)^n \beta_{2n} \pi^{2n}$ (4.4.6)

Note that the transform $g_n = T^n(f)$ can be explained by a multiple integral of order $n$ :

$$g_n(x) = \int_0^1 \int_0^1 \ldots \int_0^1 [\sum_{k=1}^{k=n} \frac{f(t_k)}{(t_k - x)\prod_{i \neq k}(t_k - t_i)} + \frac{f(x)}{\prod_{i=1}^{i=n}(x - t_i)}]dt_1 dt_2 \ldots dt_n$$

So, for $f = \varphi_0$, (4.4.6) leads to : $\lambda_{2n}^1 \pi^{2n} = \int_0^1 \int_0^1 \ldots \int_0^1 (\sum_{k=1}^{k=2n} \frac{f(t_k)}{\prod_{i \neq k}(t_k - t_i)})dt_1 dt_2 \ldots dt_{2n}$

Using basic symmetries and by change : $t_1 \leftarrow 1-t_1$ ; $t_2 \leftarrow 1-t_2$ ; ........; $t_n \leftarrow 1-t_n$ , we get

$$\int_0^1 \int_0^1 \ldots \int_0^1 (\sum_{k=1}^{k=2n} \frac{\ln(1-t_k)}{\prod_{i \neq k}(t_k - t_i)})dt_1 dt_2 \ldots dt_{2n} = -\int_0^1 \int_0^1 \ldots \int_0^1 (\sum_{k=1}^{k=2n} \frac{\ln(t_k)}{\prod_{i \neq k}(t_k - t_i)})dt_1 dt_2 \ldots dt_{2n}$$

So after simplifications:

$$\boxed{\int_0^1 \int_0^1 \ldots \int_0^1 (\sum_{k=1}^{k=2n} \frac{\ln(t_k)}{\prod_{i \neq k}(t_k - t_i)})dt_1 dt_2 \ldots dt_{2n} = \frac{\lambda_{2n}^1 \pi^{2n}}{2}}$$ (4.4.7)

## 5. Generalization of the equation.

The uniforme Lebesgue measure studied in the previous paragraph is a particular case of the following result : [See 5]

Let $\rho$ a probability density over the interval $I$ of a reducible measure associated with he reducer $\varphi(x) = \lim_{\varepsilon \to 0^+} 2\int_I \frac{(x-t)\rho(t)dt}{(x-t)^2 + \varepsilon^2}$ .

We note $T_\rho$ the operator : $f(x) \mapsto g(x) = \int_I \frac{f(t) - f(x)}{t - x}\rho(t)dt$ and $\mu(x) = \frac{\rho(x)}{\frac{\varphi^2(x)}{4} + \pi^2 \rho^2(x)}$

For any element $f$ of $L^2(I,\rho)$ that $\varphi \times f$ et $T_\rho(f)$ belong to $L^2(I,\rho)$ , and $\frac{\rho}{\mu} \times f$ is an element of $L^2(I,\mu)$, we get the formula : $\boxed{T_\rho(\varphi \times f - T_\rho(f)) = \frac{\rho}{\mu} \times f}$ (5.1)

We then define the three operators :
$a = \frac{T_\rho}{\pi}$ ; $f \xrightarrow{b} b(f) = \frac{\varphi}{2\pi} \times f$ , $f \xrightarrow{c} c(f) = \rho \times f$

With these notations, the formula (5.1) simply results in the equation :

$$\boxed{2ab = a^2 + b^2 + c^2} \quad (5.2)$$

Note that $a$ ne commute does not commute with $b$ and not usually with $c$, but $b$ and $c$ commute so obvious.

We will show by recurrence on the integer $n$ the existence of a dual system of universal coefficients $(n,k) \mapsto (u_n^k, v_n^k)$ such as for every $n \geq 1$:

$$(n+1)ab^n = a^{n+1} + nb^{n+1} + \sum_{k=1}^{k=E(\frac{n+1}{2})} (u_n^k a^{n+1-2k} + v_n^k b^{n+1-2k})c^{2k} \quad (5.3)$$

(Here $x \mapsto E(x)$ refers to the floor function).

_ For initialisation, we have : $u_1^1 = 0$ et $v_1^1 = 1$

_ Heridity takes the approach of the first paragraph, but the writing is more delicate because of the non-switching with $b, c$. Here are the details :

Multiplying left (5.3) by $a$ gives:

$$(n+1)a^2b^n = a^{n+2} + nab^{n+1} + \sum_{k=1}^{k=E(\frac{n+1}{2})} (u_n^k a^{n+2-2k} + v_n^k ab^{n+1-2k})c^{2k} \quad (5.4)$$

Multipliyng right (5.2) by $(n+1)b^n$ gives thanks the commutation of $b$ with $c$ :

$$2(n+1)ab^{n+1} = (n+1)a^2b^n + (n+1)b^{n+2} + (n+1)b^n c^2$$

Using the expression of $(n+1)a^2b^n$ in (5.4), the formula above becomes :

$$2(n+1)ab^{n+1} = a^{n+2} + nab^{n+1} + \sum_{k=1}^{k=E(\frac{n+1}{2})} (u_n^k a^{n+2-2k} + v_n^k ab^{n+1-2k})c^{2k} + (n+1)b^{n+2} + (n+1)b^n c^2$$

And by obvious groupings :

$$(n+2)ab^{n+1} = a^{n+2} + (n+1)b^{n+2} + \sum_{k=1}^{k=E(\frac{n+1}{2})} (u_n^k a^{n+2-2k} + v_n^k ab^{n+1-2k})c^{2k} + (n+1)b^n c^2 \quad (5.5)$$

For each index $k$ of the above sum, we have by induction hypothesis :

$$ab^{n+1-2k} = \frac{1}{n+2-2k}[a^{n+2-2k} + (n+1-2k)b^{n+2-2k} + \sum_{j=1}^{j=E(\frac{n+2-2k}{2})} (u_{n+1-2k}^j a^{n+2-2k-2j} + v_{n+1-2k}^j b^{n+2-2k-2j})c^{2j}]$$

Multiplying right by $c^{2k}$ and withe the change $q = k + j$, we have :

$$\sum_{j=1}^{j=E(\frac{n+2-2k}{2})}(u_{n+1-2k}^j a^{n+2-2k-2j} + v_{n+1-2k}^j b^{n+2-2k-2j})c^{2j}]c^{2k} = \sum_{q=k+1}^{q=E(\frac{n+2}{2})}(u_{n+1-2k}^{q-k} a^{n+2-2q} + v_{n+1-2q}^{q-k} b^{n+2-2q})c^{2q}$$

By substituting in (5.5) the expressions of $ab^{n+2-2k}c^{2k}$ obtained above, we can write:

$$(n+2)ab^{n+1} = a^{n+2} + (n+1)b^{n+2} + (n+1)b^n c^2 + S_1 + S_2 + S_3 \text{ with :}$$

$$S_1 = \sum_{k=1}^{k=E(\frac{n+1}{2})} u_n^k a^{n+2-2k} c^{2k} \quad ; \quad S_2 = \sum_{k=1}^{k=E(\frac{n+1}{2})} \frac{v_n^k}{n+2-2k}[a^{n+2-2k} + (n+1-2k)b^{n+2-2k}]c^{2k}$$

$$S_3 = \sum_{q=2}^{q=E(\frac{n+2}{2})} \left( \sum_{k=1}^{k=q-1} \frac{v_n^k}{n+2-2k}[u_{n+1-2k}^{q-k} a^{n+2-2q} + v_{n+1-2k}^{q-k} b^{n+2-2q}] \right) c^{2q}$$

This is the expression searched for the next rank $n+1$:

$$(n+2)ab^{n+1} = a^{n+2} + (n+1)b^{n+2} + \sum_{k=1}^{k=E(\frac{n+2}{2})}[u_{n+1}^k a^{n+2-2k} + v_{n+1}^k b^{n+2-2k}]c^{2k}$$

The new coefficients are given by the following relations:

- For $v_{n+1}^q$:

$$v_{n+1}^1 = n+1 + \frac{(n-1)}{n} v_n^1$$

If $2 \leq q \leq E(\frac{n+1}{2})$ $\quad v_{n+1}^q = \sum_{k=1}^{k=q-1} \frac{v_n^k v_{n+1-2k}^{q-k}}{n+2-2k} + \frac{(n+1-2q)}{n+2-2q} v_n^q$

If $n$ is even: $v_{n+1}^{\frac{n+2}{2}} = \sum_{k=1}^{k=E(\frac{n+1}{2})} \frac{v_n^k v_{n+1-2k}^{\frac{n+2}{2}-k}}{n+2-2k}$

- For $u_{n+1}^q$:

$$u_{n+1}^1 = u_n^1 + \frac{v_n^1}{n}$$

If $2 \leq q \leq E(\frac{n+1}{2})$ $\quad u_{n+1}^q = \sum_{k=1}^{k=q-1} \frac{v_n^k u_{n+1-2k}^{q-k}}{n+2-2k} + \frac{v_n^q}{n+2-2q} + u_n^q$

If $n$ is even: $u_{n+1}^{\frac{n+2}{2}} = \sum_{k=1}^{k=E(\frac{n+1}{2})} \frac{v_n^k u_{n+1-2k}^{\frac{n+2}{2}-k}}{n+2-2k}$

In the particular case of the Lebesgue measure, we have $c=1$, and so we get :

$$ab^n = \frac{1}{n+1}[a^{n+1} + \sum_{k=1}^{k=E(\frac{n+1}{2})} u_n^k a^{n+1-2k} + nb^{n+1} + \sum_{k=1}^{k=E(\frac{n+1}{2})} v_n^k b^{n+1-2k}]$$

Or in this situation we have obtained : $ab^n = A_n(a) + C_n(b)$

So we deduce :

$$A_n(X) = \sum_{k=1}^{k=n+1} \alpha_n^k X^k = \frac{1}{n+1}(\sum_{q=0}^{q=E(\frac{n+1}{2})} u_n^q X^{n+1-2q} + X^{n+1})$$

$$B_n(X) = \sum_{k=0}^{k=n+1} \lambda_n^k X^k = \frac{1}{n+1}(\sum_{q=1}^{q=E(\frac{n+1}{2})} v_n^q X^{n+1-2q} + nX^{n+1})$$

And finally, $\boxed{\begin{cases} u_n^k = (n+1)\alpha_n^{n+1-2k} \\ v_n^k = (n+1)\lambda_n^{n+1-2k} \end{cases}}$ (5.6)

Application to the calculation of integrals of type $\int_0^1 \varphi^n(t)\rho(t)dt$

The formula (5.3) gets here :

$$(n+1)\frac{T_\rho}{\pi}(\frac{\varphi^n}{(2\pi)^n} \times f) = (\frac{T_\rho}{\pi})^{n+1}(f) + n\left(\frac{\varphi}{2\pi}\right)^{n+1} f + \sum_{k=1}^{k=E(\frac{n+1}{2})} [u_n^k(\frac{T_\rho}{\pi})^{n+1-2k}(\rho^{2k}f) + v_n^k(\frac{\varphi}{2\pi})^{n+1-2k}\rho^{2k}f]$$

With $u_n^0 = 1$ ; $v_n^0 = n$ and $\varphi_0 = \frac{\varphi}{2}$, this simplifies to:

$$\boxed{(n+1)T_\rho(\varphi_0^n f) = \sum_{k=0}^{k=E(\frac{n+1}{2})} \pi^{2k}[u_n^k T_\rho^{n+1-2k}(\rho^{2k}f) + v_n^k \rho^{2k}\varphi_0^{n+1-2k} f]} \quad (5.7)$$

For $f = 1$ we get : $(n+1)T_\rho(\varphi_0^n) = \sum_{k=0}^{k=E(\frac{n+1}{2})} \pi^{2k}[u_n^k T_\rho^{n+1-2k}(\rho^{2k}) + v_n^k \rho^{2k}\varphi_0^{n+1-2k}]$

According to the method used in paragraph (4) we can evaluate $(n+1)T_\rho(\varphi_0^n x)$ by two different ways :

➢ By noting $f = x$ in (5.7).

➢ By using : $T_\rho(x\varphi_0^n) = xT_\rho(\varphi_0^n) + \int_I \varphi_0^n(t)\rho(t)dt$

_ The first method leads to :

$$(n+1)T_\rho(\varphi_0^n x) = \sum_{k=0}^{k=E(\frac{n+1}{2})} \pi^{2k}[u_n^k T_\rho^{n+1-2k}(\rho^{2k} x) + v_n^k \rho^{2k} \varphi_0^{n+1-2k} x]$$

Or, for every $q$ integer: $T_\rho^q(x\rho^{2k}) = xT_\rho^q(\rho^{2k}) + \int_I T_\rho^{q-1}(\rho^{2k})(t)\rho(t)dt$ , so we get :

$$(n+1)T_\rho(\varphi_0^n x) = \sum_{k=0}^{k=E(\frac{n+1}{2})} \pi^{2k}[u_n^k x T_\rho^{n+1-2k}(\rho^{2k}) + v_n^k \rho^{2k} \varphi_0^{n+1-2k} x] + S$$

with $S = \sum_{k=0}^{k=E(\frac{n+1}{2})} \pi^{2k} u_n^k \int_I T_\rho^{n-2k}(\rho^{2k})(t)\rho(t)dt$

_ The secon method gives :

$$(n+1)T_\rho(\varphi_0^n x) = x \sum_{k=0}^{k=E(\frac{n+1}{2})} \pi^{2k}[u_n^k T_\rho^{n+1-2k}(\rho^{2k}) + v_n^k \rho^{2k} \varphi_0^{n+1-2k}] + (n+1)\int_I \varphi_0^n(t)\rho(t)dt$$

By comparing the two results we finally obtain :

$$\boxed{(n+1)\int_I \varphi_0^n(t)\rho(t)dt = \sum_{k=0}^{k=E(\frac{n+1}{2})} \pi^{2k} u_n^k \int_I T_\rho^{n-2k}(\rho^{2k})(t)\rho(t)dt} \quad (5.8)$$

References.

[1] Steven Roman et Gian-Carlo Rota, *"The Umbral Calculus"*, *Advances in Mathematics*, volume 27, pages 95 - 188, (1978). (1978).

[2] Milton Abramowitz and Irene A. Stegun, eds. *Handbook of Mathematical Functions* with Formulas, Graphs, *and Mathematical Tables*, (1972) Dover, New York.

[3] William C. Brenke, *On generating functions of polynomial systems*, « American Mathematical Monthly », (1945) **52** pp. 297-301.

[4] Harris, J. W. and Stocker, H. "Secant and Cosecant."  *Handbook of Mathematics and Computational Science.* New York: Springer-Verlag, 1998.

[5] Roland Groux.  arXiv :1104.3218v1 [math.CA]